\documentclass[12pt]{amsart}

\usepackage{amsmath,xspace,amssymb,mathrsfs}
\usepackage{color}

\input xy
\xyoption{all}
\xyoption{2cell}
\UseAllTwocells
\CompileMatrices

\newcommand{\Spec}{\operatorname{Spec}}
\renewcommand{\phi}{\varphi}

\newcommand{\Clop}{\operatorname{Clop}}

\newtheorem{proposition}{Proposition}[section]
\newtheorem{lemma}[proposition]{Lemma}

\newtheorem{corollary}[proposition]{Corollary}
\newtheorem{theorem}[proposition]{Theorem}

\theoremstyle{definition}
\newtheorem{definition}[proposition]{Definition}

\begin{document}

\title[Ring-theoretic approaches to topology]{Ring-theoretic approaches to point-set topology}
\author[A. Tarizadeh]{Abolfazl Tarizadeh}
\address{Department of Mathematics, Faculty of Basic Sciences, University of Maragheh \\
P. O. Box 55136-553, Maragheh, Iran.}
\email{ebulfez1978@gmail.com}

\date{}
\footnotetext{2010 Mathematics Subject Classification: 54A20, 14A05, 13B10, 54B10, 54B99.
\\ Keywords: Power set ring; Zariski convergent; Tychonoff theorem; Alexander subbase theorem; Profinite space.\\}

\begin{abstract} In this paper, it is shown that a topological space $X$ is compact iff every maximal ideal of the power set ring $\mathcal{P}(X)$ converges to exactly one point of $X$. Then as an application, simple and ring-theoretic proofs are provided for the Tychonoff theorem and Alexander subbase theorem. As another result in this spirit, a ring-theoretic proof is given to the fact that a topological space is a profinite space iff it is compact and totally disconnected. \\
\end{abstract}

\maketitle

\section{Introduction}

Some of the mathematicians consider the Tychonoff theorem as the single most important result in general topology, (other mathematicians allow it to share this honor with Urysohn's lemma). Tychonoff theorem is a fundamental ingredient in proving various important results in the areas of topology, geometry, algebra and mathematical analysis. This result has been investigated in the literature over the years from various point of view, see e.g. \cite{Celmentino}- \cite{Wright}. Alexander subbase theorem and characterization of profinite spaces are another fundamental results in general topology which have significant applications in mathematics. \\
In this paper, an algebraic characterization for the quasi-compactness of a topological space is given, see Theorem \ref{Theorem I}. Using this, then the Tychonoff theorem and Alexander subbase theorem are easily deduced. An algebraic characterization for the Hausdorffness of a topological space is also given, see Theorem \ref{Theorem III}. As an immediate consequence of Theorems \ref{Theorem I} and \ref{Theorem III}, we obtain that a topological space $X$ is compact if and only if every maximal ideal of the power set ring $\mathcal{P}(X)$ converges to exactly one point of $X$. Finally, we give an algebraic proof to the fact that a topological space is a profinite space if and only if it is compact and totally disconnected, see Theorem \ref{Theorem II}. Most of the proofs are greatly based on the Zariski convergent notion and on the significant using of the power set ring. \\

\section{Preliminaries}

If $X$ is a set, then its power set $\mathcal{P}(X)$ together with the symmetric difference $A+B=(A\cup B)\setminus (A\cap B)$ as the addition and the intersection $A.B=A\cap B$ as the multiplication forms a commutative ring whose zero and unit are respectively the empty set and the whole set $X$. The ring $\mathcal{P}(X)$ is called the \emph{power set ring} of $X$. If $f:X\rightarrow Y$ is a function, then the map $\mathcal{P}(f):\mathcal{P}(Y)\rightarrow\mathcal{P}(X)$ defined by $A\rightsquigarrow f^{-1}(A)$ is a morphism of rings. In fact, the assignments $X\rightsquigarrow\mathcal{P}(X)$ and $f\rightsquigarrow\mathcal{P}(f)$ form a faithful contravariant functor from the category of sets to the category of commutative rings. We call it the \emph{power set functor}. \\
A ring is called a Boolean ring if every element is idempotent. It is easy to see that every Boolean ring is a commutative ring, and in a Boolean ring every prime ideal is a maximal ideal. The power set ring $\mathcal{P}(X)$ is a typical example of Boolean rings. \\
If $f_{1},...,f_{n}$ are a finite number of idempotent elements of a commutative ring $R$, then $(f_{1},...,f_{n})$ is a principal ideal of $R$ and generated by an idempotent, since $(f_{1},f_{2})=(f_{1}+f_{2}-f_{1}f_{2})$. \\
If $A\in\mathcal{P}(X)$ then $\mathcal{P}(A)$ is an ideal of $\mathcal{P}(X)$. In fact, $\mathcal{P}(A)=(A)$ is a principal ideal. If $A_{1},...,A_{n}$ are a finite number of elements of $\mathcal{P}(X)$ then $(A_{1},...,A_{n})=\mathcal{P}(\bigcup\limits_{i=1}^{n}A_{i})$. \\
If $x\in X$, then clearly $\mathfrak{m}_{x}:=\mathcal{P}(X\setminus\{x\})$ is a maximal ideal of $\mathcal{P}(X)$. It is also easy to see that $X$ is a finite set iff every maximal ideal of $\mathcal{P}(X)$ is of the form $\mathfrak{m}_{x}$. \\
If $\phi:A\rightarrow B$ is a morphism of rings, then the induced map $\Spec(B)\rightarrow\Spec(A)$ is given by $\mathfrak{p}\rightsquigarrow\phi^{-1}(\mathfrak{p})$. \\
If $f$ is a member of a ring $R$, then $D(f)=\{\mathfrak{p}\in\Spec(R): f\notin\mathfrak{p}\}$. \\

\begin{definition} Let $X$ be a topological space, $x\in X$ and $M$ a maximal ideal of $\mathcal{P}(X)$. Then we say that $M$ is \emph{convergent} (or, \emph{Zariski convergent}) to the point $x$ if $U$ is an open of $X$ containing $x$ then $M\in D(U)$. \\
\end{definition}

Let $\phi:X\rightarrow Y$ be a continuous map of topological spaces. If a maximal ideal $M$ of $\mathcal{P}(X)$ converges to some point $x\in X$, then clearly the maximal ideal $\mathcal{P}(\phi)^{-1}(M)$ of $\mathcal{P}(Y)$ is convergent to $\phi(x)$. \\Pl\"{u}cker
A quasi-compact and Hausdorff topological space is called a compact space. \\

\begin{definition}\label{Definition I} Let $(X_{i},\phi_{i,j})$ be a projective (inverse) system of finite discrete spaces over a poset $(I,<)$. The projective (inverse) limit $X=\lim\limits_{\overleftarrow{i\in I}}X_{i}$ with the induced product topology, as a subset of $\prod\limits_{i\in I}X_{i}$, is called a \emph{profinite space}. \\
\end{definition}

\section{Main results}

The following is one of the main results of this paper which builds a bridge between topology and commutative ring theory. \\

\begin{theorem}\label{Theorem I} A topological space $X$ is quasi-compact if and only if every maximal ideal of $\mathcal{P}(X)$ converges to a point of $X$. \\
\end{theorem}

{\bf Proof.} Let $X$ be quasi-compact and let $M$ be a maximal ideal of $\mathcal{P}(X)$. The collection of closures $\overline{A}$ with $A\in S:=\mathcal{P}(X)\setminus M$ has the finite intersection property. Therefore $\bigcap\limits_{A\in S}\overline{A}\neq\emptyset$, since $X$ is quasi-compact. Thus we may choose some point $x$ in the intersection. If $U$ is an open of $X$ containing $x$, then it will be enough to show that $U\notin M$. If $U\in M$ then
$U^{c}=X\setminus U\notin M$ and so $x\in\overline{U^{c}}=U^{c}$ which is a contradiction. Conversely, if $\{U_{i}: i\in I\}$ is an open covering of $X$ then we claim that the ideal $(U_{i}: i\in I)$ is the whole ring $\mathcal{P}(X)$. If not, then there exists a maximal ideal $M$ of $\mathcal{P}(X)$ containing this ideal. By the hypothesis, $M$ converges to a point $x\in X$. Clearly $x\in U_{k}$ for some $k$. It follows that $U_{k}\notin M$ which is a contradiction. This establishes the claim. Thus there exists a finite subset $J\subseteq I$ such that $\mathcal{P}(X)=(U_{i}: i\in J)=\mathcal{P}(\bigcup\limits_{i\in J}U_{i})$. This yields that $X=\bigcup\limits_{i\in J}U_{i}$. $\Box$ \\

\begin{corollary}$($Tychonoff Theorem$)$ Let $(X_{i})_{i\in I}$ be a family of quasi-compact topological spaces. Then $X=\prod\limits_{i\in I}X_{i}$ with the product topology is quasi-compact. \\
\end{corollary}

{\bf Proof.} Let $M$ be a maximal ideal of $\mathcal{P}(X)$. Setting $M_{i}:=\mathcal{P}(\pi_{i})^{-1}(M)$ where $\pi_{i}:X\rightarrow X_{i}$ is the canonical projection map. For each $i\in I$ then by Theorem \ref{Theorem I}, the maximal ideal $M_{i}$ is convergent to a point $x_{i}\in X_{i}$. To prove the assertion, by Theorem \ref{Theorem I}, it suffices to show that $M$ is convergent to the point $x=(x_{i})$. Let $U$ be an open of $X$ containing $x$. Then there exists a basis open $V=\prod\limits_{i\in I}V_{i}$ of $X$ such that $x\in V\subseteq U$ where each $V_{i}$ is an open of $X_{i}$ and $V_{i}=X_{i}$ for all but a finite number of indices $i$. Let $J$ be the set of all $i\in I$ such that $V_{i}\neq X_{i}$. Then $V=\bigcap\limits_{i\in J}\pi^{-1}(V_{i})$. Clearly $V_{i}\notin M_{i}$ and so $\pi^{-1}_{i}(V_{i})\notin M$ for all $i$. Thus $V\notin M$ since $J$ is a finite set. Therefore $U\notin M$. $\Box$ \\

\begin{corollary}$($Alexander Subbase Theorem$)$ Let $X$ be a topological space and let $\mathscr{D}$ be a subbase of $X$ such that every covering of $X$ by elements of $\mathscr{D}$ has a finite refinement. Then $X$ is quasi-compact. \\
\end{corollary}

{\bf Proof.} Let $M$ be a maximal ideal of $\mathcal{P}(X)$. By Theorem \ref{Theorem I}, it will be enough to show that $M$ is convergent to a point of $X$. By the hypothesis, $X\neq\bigcup\limits_{D\in M\cap\mathscr{D}}D$ since otherwise we may find a finite number $D_{1},...,D_{n}\in M\cap\mathscr{D}$ such that $X=\bigcup\limits_{i=1}^{n}D_{i}$ and so $\mathcal{P}(X)=(D_{1},...,D_{n})\subseteq M$ which is a contradiction. Hence, we may choose some $x\in X$ such that $x\notin\bigcup\limits_{D\in M\cap\mathscr{D}}D$. If $U$ is an open of $X$ containing $x$ then there exists a finite number $D'_{1},...,D'_{s}\in\mathscr{D}$ such that $x\in\bigcap\limits_{k=1}^{s}D'_{k}\subseteq U$. But $\bigcap\limits_{k=1}^{s}D'_{k}\notin M$. Therefore $U\notin M$. $\Box$ \\

\begin{theorem}\label{Theorem III} A topological space $X$ is Hausdorff if and only if every maximal ideal of $\mathcal{P}(X)$ converges to at most one point of $X$. \\
\end{theorem}

{\bf Proof.} Assume that $X$ is Hausdorff. Suppose there exists a maximal ideal $M$ of $\mathcal{P}(X)$ which is convergent to the distinct points $x$ and $y$ of $X$. Then we may find opens $U$ and $V$ of $X$ with $x\in U$, $y\in V$ and $U\cap V=\emptyset$. It follows that either $U\in M$ or $V\in M$, which is a contradiction. Conversely, let $x$ and $y$ be two distinct points of $X$. Let $I$ be the ideal of $\mathcal{P}(X)$ generated by all $U^{c}$ where $U$ is an open of $X$ containing $x$. Similarly, let $J$ be the ideal of $\mathcal{P}(X)$ generated by all $V^{c}$ where $V$ is an open of $X$ containing $y$. We claim that $I+J=\mathcal{P}(X)$. If not, then there exists a maximal ideal $M$ of  $\mathcal{P}(X)$ such that $I+J\subseteq M$. It follows that $M$ is convergent to the both points $x$ and $y$, which is a contradiction. This establishes the claim. Thus we may find a finite number $U_{1},...,U_{n}$ of open neighborhoods of $x$ with $n\geq1$ and a finite number $V_{1},...,V_{m}$ of open neighborhoods of $y$ with $m\geq1$ such that $\mathcal{P}(X)=(U^{c}_{1},...,U^{c}_{n})+(V^{c}_{1},...,V^{c}_{m})$. This yields that $(\bigcap\limits_{i=1}^{n}U_{i})\cap
(\bigcap\limits_{j=1}^{m}V_{j})=\emptyset$. $\Box$ \\

\begin{corollary} A topological space $X$ is compact if and only if every maximal ideal of $\mathcal{P}(X)$ converges to exactly one point of $X$. \\
\end{corollary}

{\bf Proof.} It is an immediate consequence of Theorems \ref{Theorem I} and \ref{Theorem III}. $\Box$ \\

Let $R$ be a Boolean ring. Then $\{0,1\}$ is a subring of $R$. If $A$ is a subring of $R$ and $f\in R$ then $A[f]=\{a+bf: a,b\in A\}$. In particular, if $A$ is a finite subring of $R$ and $f_{1},...,f_{n}\in R$ then $A[f_{1},...,f_{n}]$ is also a finite subring of $R$. Let $\{R_{i}: i\in I\}$ be the set of all finite subrings of $R$.
We define $j<i$ if $R_{j}$ is a proper subset of $R_{i}$. Then the poset $(I,<)$ is directed, because if $A$ and $B$ are finite subrings of $R$ then we observed in the above that $A[B]$ is a finite subring of $R$ and $A,B\subseteq A[B]$. \\

\begin{lemma}\label{Lemma I} Let $R$ be a Boolean ring. Then $\Spec(R)$ is a profinite space. \\
\end{lemma}

{\bf Proof.} Let $\{R_{i}: i\in I\}$ be the set of all finite subrings of $R$.
Then the $\Spec(R_{i})$ together with the $\phi_{i,j}:\Spec(R_{i})\rightarrow\Spec(R_{j})$ induced by the inclusions $R_{j}\subseteq R_{i}$, as the transition morphisms, form a projective system of finite discrete spaces over the poset $(I,<)$. We show that $\Spec(R)$ together with the canonical maps $p_{i}:\Spec(R)\rightarrow\Spec(R_{i})$, induced by the inclusions $R_{i}\subseteq R$, is the projective limit of the above system. By the universal property of the projective limits, there exists a (unique) continuous map $\phi:\Spec(R)\rightarrow X=\lim\limits_{\overleftarrow{i\in I}}\Spec(R_{i})$ such that $p_{i}=\pi_{i}\circ\phi$ for all $i$, where each $\pi_{i}:X\rightarrow\Spec(R_{i})$ is the canonical projection. Therefore $\phi(M)=(M\cap R_{i})$. The map $\phi$ is clearly a closed map because $\Spec(R)$ is quasi-compact and $X$ is Hausdorff. It remains to show that it is bijective. Suppose $\phi(M)=\phi(N)$. If $f\in R$ then $A=\{0,1,f,1+f\}$ is a subring of $R$ and so $M\cap A=N\cap A$. Thus
$M= N$. Finally, take $(M_{i})\in X$ where each $M_{i}$ is a maximal ideal of $R_{i}$. Let $M$ be the ideal of $R$ generated by the subset $\bigcup\limits_{i\in I}M_{i}\subseteq R$. We show that $M$ is a maximal ideal of $R$. Clearly $M$ is a proper ideal of $R$. If not, then there exists a finite subset $J\subseteq I$ such that the ideal generated by the subset $\bigcup\limits_{i\in J}M_{i}$ is the whole ring $R$.
But we may find some $k\in I$ such that $i\leq k$ for all $i\in J$,
since $(I,<)$ is directed. We have $\bigcup\limits_{i\in J}M_{i}\subseteq M_{k}$ because $M_{i}=R_{i}\cap M_{k}$ for all $i\in J$. This yields that $M_{k}=R_{k}$ which is a contradiction. If $f,g\in R$ such that $fg\in M$ then similarly above there exists some $k\in I$ such that $f,g \in R_{k}$ and $fg\in M_{k}$. Thus either $f\in M_{k}$ or $g\in M_{k}$. Therefore $M$ is a prime ideal of $R$. Clearly $M_{i}\subseteq M\cap R_{i}$ and so $M_{i}=M\cap R_{i}$ for all $i\in I$. $\Box$ \\

\begin{lemma}\label{Lemma II} Let $(X_{i})$ be a family of Hausdorff and totally disconnected spaces. If $\prod\limits_{i}X_{i}$ is equipped with the product topology, then every subspace $X$ is Hausdorff and totally disconnected. \\
\end{lemma}

{\bf Proof.} It is clearly Hausdorff. Let $C$ be a connected subset of $X$. If $(x_{i})$ and $(y_{i})$ are two distinct points of $C$ then there exists some $k$ such that $x_{k}\neq y_{k}$. Clearly $x_{k},y_{k}\in\pi_{k}(C)$ and $\pi_{k}(C)$ is a connected subset of $X_{k}$ where $\pi_{k}:X\rightarrow X_{k}$ is the canonical projection. So $\pi_{k}(C)$ is a single-point set. But this is a contradiction and we win. $\Box$ \\

\begin{theorem}\label{Theorem II} A topological space is a profinite space if and only if it is compact and totally disconnected. \\
\end{theorem}

{\bf Proof.} Let $X$ be a profinite space. By Lemma \ref{Lemma II}, it is Hausdorff and totally disconnected. To see the quasi-compactness we use Theorem \ref{Theorem I}. So let $M$ be a maximal ideal of $\mathcal{P}(X)$. By taking into account the notations of Definition \ref{Definition I}, setting  $M_{i}:=\mathcal{P}(\pi_{i})^{-1}(M)$ where each $\pi_{i}:X\rightarrow X_{i}$ is the canonical projection. But each $X_{i}$ is a finite set and so  $M_{i}=\mathcal{P}(X_{i}\setminus\{x_{i}\})$ for some $x_{i}\in X_{i}$. We prove that $M$ is convergent to $x=(x_{i})$. First we have to show that $x\in X$. If $j\leq i$ then $\pi_{j}=\phi_{i,j}\circ\pi_{i}$. It follows that $M_{j}=\mathcal{P}(\phi_{i,j})^{-1}(M_{i})$. This yields that $M_{j}$ is convergent to $\phi_{i,j}(x_{i})$ because $M_{i}$ is obviously convergent to $x_{i}$. Therefore $\phi_{i,j}(x_{i})=x_{j}$. Hence, $x\in X$. Now let $U$ be an open of $X$ containing $x$. There exists a basis open $V=\prod\limits_{i\in I}V_{i}$ in the product topology such that $x\in X\cap V\subseteq U$. Thus there exists a finite subset $J\subseteq I$ such that $X\cap V=\bigcap\limits_{i\in J}\pi^{-1}(V_{i})$ and $\pi^{-1}(V_{i})\notin M$ for all $i$. Therefore $X\cap V$ and so $U$ are not in $M$. Hence, $M$ is convergent to $x$. Conversely, let $X$ be a compact totally disconnected space. Let $R=\Clop(X)$ be the set of all clopen (both open and closed) subsets of $X$. Clearly $R$ is a subring of $\mathcal{P}(X)$. It is well known that for each point $x$ in a compact space $X$, then the connected component of $X$ containing $x$ is the intersection of all $A\in R$ such that $x\in A$. Using this fact, then it can be shown that the map $X\rightarrow\Spec(R)$ given by $x\rightsquigarrow\mathfrak{m}_{x}\cap R$ is a homeomorphism where $\mathfrak{m}_{x}=\mathcal{P}(X\setminus\{x\})$. Therefore by Lemma \ref{Lemma I}, $X$ is a profinite space. $\Box$ \\

\end{document}